\documentclass[12pt]{amsart}
\usepackage{amsmath,amsfonts,amsthm,amscd}
\newtheorem{theorem}{Theorem}[section]
\newtheorem{claim}[theorem]{Claim}
\newtheorem{proposition}[theorem]{Proposition}
\newtheorem{lemma}[theorem]{Lemma}
\newtheorem{definition}[theorem]{Definition}
\newtheorem{corollary}[theorem]{Corollary}
\theoremstyle{remark}
\newtheorem*{remark}{Remark}
\numberwithin{equation}{section}
\begin{document}
\bibliographystyle{amsalpha}
\title{Volume growth, curvature decay, and critical metrics} 
\author{Gang Tian}
\address{Gang Tian\\ Fine Hall, Princeton University, 
Princeton, NJ 08544}
\email{tian@math.princeton.edu}
\thanks{The research of the first author was 
partially supported by NSF Grant DMS-0302744.}
\author{Jeff Viaclovsky}
\address{Jeff Viaclovsky, Department of Mathematics, 
MIT, Cambridge, MA 02139}
\address{Department of Mathematics, University of Wisconsin, Madison, WI, 53706}
\email{jeffv@math.wisc.edu}
\thanks{The research of the second author was partially 
supported by NSF Grant DMS-0503506.}
\date{December 17, 2006}
\begin{abstract}
We make some improvements to our previous results 
in \cite{TV} and \cite{TV2}. First, we prove a version 
of our volume growth theorem which does not require 
any assumption on the first Betti number. Second, we show that 
our local regularity theorem only requires 
a lower volume growth assumption, not a full 
Sobolev constant bound. As an application of these 
results, we can weaken the assumptions of several of our theorems in
\cite{TV} and \cite{TV2}.  
\end{abstract}
\maketitle
%%%%%%%%%%%%%%%%%%%%%%%%%%%%%%%%%%%%%%%%%%%%%%%%%%%%%%
\section{Introduction}

 Riemannian spaces with quadratic curvature decay have been widely 
studied in the literature, see for example \cite{Abresch1},
\cite{Abresch2} , \cite{AbreschGromoll}, \cite{BKN}, \cite{Kasue1}, 
\cite{Kasue2}, \cite{GreenePetersen}, \cite{Gromovdiam}, \cite{Zhu}. All
of these works assume that the curvature 
decay is strictly better than quadratic in the sense that
\begin{align}
|Rm| = O(r^{-(2+ \delta)}), \mbox{ as } r \rightarrow \infty,
\end{align}
for some $\delta > 0$, where $r(x) = d(p,x)$ is 
the distance to a basepoint $p$, or the weaker assumption that
\begin{align}
\label{kdecay}
Rm \geq - \frac{k(r)}{r^2},
\end{align}
(meaning all of the sectional curvatures are bounded below 
accordingly), and the function $k(r)$ satisfies 
\begin{align}
\label{inta}
\int_1^{\infty} \frac{ k(r)}{r} dr < \infty.
\end{align}
Spaces satifying such a curvature decay condition  
are said to have {\em asymptotically nonnegative curvature}. 
We remark that by standard comparsion theory, (\ref{kdecay}) 
and (\ref{inta}) imply an upper volume growth estimate, 
\begin{align}
Vol(B(p,r)) \leq V_0 r^n,
\end{align}
for some constant $V_0$. Moreover, only a lower bound on the Ricci 
curvature is needed for this upper volume growth estimate \cite{Zhu}.

 In our investigation of critical metrics in \cite{TV}, 
\cite{TV2}, and in the work of \cite{Andersonc},  
spaces arise with curvature decay as
in (\ref{kdecay}), but the function $k(r)$ only 
satisfies $k(r) \rightarrow 0$ as $r \rightarrow \infty$,
and standard comparsion arguments do not apply.
In \cite{TV} we proved an upper volume growth estimate in
this case, but our proof required finiteness of the first Betti number 
to rule out the presence of so-called ``bad'' annuli. In this paper, 
we show that adding the condition (\ref{csbnd}) below, eliminates this pathology.
For $M$ non-compact, $C_S$ is defined to be the best constant so that 
\begin{align}
\label{mainSob2}
\Vert f \Vert_{L^{\frac{2n}{n-2}}} \leq C_S   
\Vert \nabla f \Vert_{L^2},
\end{align}
for all $f \in C^{0,1}(M)$ with compact support. 
Let $Ric_-$ denote the negative part of the Ricci tensor. 
\begin{theorem}
\label{bigthm_j}
Let $(M,g)$ be a complete, noncompact, $n$-dimensional 
Riemannian manifold with base point $p$.
Assume that 
\begin{align}
\label{cond4_i}
C_S < \infty,
\end{align}
and that 
\begin{align}
\label{decay1_i}
\underset{S(r)}{sup} \ |Rm_g| &= o(r^{-2}),
\end{align}
 as $r \rightarrow \infty$,
where $S(r)$ denotes the sphere of radius $r$ centered at $p$. 
If 
\begin{align}
\label{csbnd}
&\int_M | Ric_-|^{\frac{n}{2}} dV_g < \Lambda,
\end{align}
for some constant $\Lambda \in \mathbb{R}$, then $(M,g)$ has finitely many ends, 
and there exists a constant $C_2$ (depending on $g$) so that 
\begin{align}
\label{vga_i}
Vol(B(p,r)) \leq C_2 r^n.
\end{align}
Furthermore, each end is ALE of order $0$.
\end{theorem}
We have another generalization of our previous results. 
We consider any system of the type
\begin{align}
\label{generaleqn}
\Delta Ric = Rm * Ric,
\end{align}
where the right hand side is shorthand notation 
for a linear combinating of terms of the form 
$A_{ijkl}B_{jl}$, where $A_{ijkl}$ depends on the 
full curvature tensor, and $B_{jl}$ depends only on
the Ricci tensor. 
We call any metric satisfying a system of the 
form (\ref{generaleqn}) a {\em{critical metric}}.
In \cite[Theorem 3.1]{TV}, we proved an $\epsilon$--regularity theorem 
which depended on the Sobolev constant. Here we relax this 
condition and require only a lower volume growth 
assumption:
\begin{theorem}
\label{higherlocalregthm}
Assume that (\ref{generaleqn}) is satisfied, 
let $r < diam(M)/2$, and $B(p,r)$ be a geodesic 
ball around the point $p$, and $k \geq 0$. 
If there exists a constant $V_0 > 0$ so 
that 
\begin{align*}
Vol( B(q,s)) \geq V_0 s^4
\end{align*}
for all $q \in B(p, r/2),$ and $s \leq r/2$,
then there exist  
constants $\epsilon_0, C_k$  (depending upon $V_0$) so that if 
\begin{align*}
\Vert Rm \Vert_{L^2(B(p,r))} = 
\left\{ \int_{B(p,r)} |Rm|^2 dV_g \right\}^{1/2} \leq \epsilon_0,
\end{align*}
then 
\begin{align*}
\underset{B(p, r/2)}{sup}| \nabla^k Rm| \leq
\frac{C_k}{r^{2+k}} \left\{ \int_{B(p,r)} |Rm|^2 dV_g \right\}^{1/2}
\leq \frac{C_k \epsilon_0}{r^{2+k}}. 
\end{align*}
\end{theorem}

 A consequence of these results is that we can (i) 
remove the Betti number assumption from our volume growth 
theorem from \cite[Theorem 1.2]{TV2}, or (ii) We can  
relax the Sobolev constant assumption to only an assumption 
on lower volume growth. 

Recall that a metric $g$ is called anti-self-dual if 
$W^+_g \equiv 0$, and self-dual if $W^-_g \equiv 0$. 
As in \cite{TV2}, we specialize to the class of \\
\begin{tabular}{ll}\vspace{-3mm}
\\
(a) & self-dual or anti-self-dual metrics with constant scalar curvature,\\

(b) & metrics with harmonic curvature ($\delta Rm \equiv 0$),\\

(c) & K\"ahler metrics with constant scalar curvature.\\ 
\end{tabular}
\\
\\
We have the following notion 
of local Sobolev constant. 
For $p \in M$, and $r > 0$, we define $C_{S}(p, r)$ to be the best 
constant such that  
\begin{align}
\label{mainSob3}
\Vert f \Vert_{L^{4}} \leq C_{S}(p,r)   
\Vert \nabla f \Vert_{L^2},
\end{align}
for all $f \in C^{0,1}$ with compact support 
in $B(p, r)$, and define 
\begin{align}
C_{S}(r) = \sup_{p \in M} C_S(p,r).
\end{align}
Let $b_1(M)$ denote the first 
Betti number of $M$. 
\begin{theorem} 
\label{orbthm2}
Let $(M, g)$ be a metric of class (a), (b), or (c) on a 
smooth, complete four-dimensional manifold $M$ satisfying
\begin{align} 
\label{l2a}
\int_{M} |Rm_{g}|^2 dV_g \leq \Lambda,
\end{align}
for some constant $\Lambda$. 

Assume that 
\begin{align} 
Vol( B(q,s)) &\geq V_0 s^4, \mbox{ for all } q \in M,  
\mbox{ and } s \leq diam(M)/2, \\
b_1(M) &< B_1,
\end{align}
where $V_0, B_1$ are constants.
Then there exists a constant 
$V_1$, depending only upon $V_0, \Lambda, B_1$, 
such that $Vol (B(p,r)) \leq V_1 \cdot r^4$, 
for all $p \in M$ and $r > 0$.  

Assume instead that  
\begin{align} 
& C_{S}(r)  < C_1, \mbox{ for } r < diam(M)/2,
\end{align}
where $C_1$ is a constant.
Then there exists a constant 
$V_2$, depending only upon $C_1, \Lambda$, 
such that $Vol (B(p,r)) \leq V_2 \cdot r^4$, 
for all $p \in M$ and $r > 0$.  
\end{theorem}

We have the following improvement 
of our compactness theorem \cite[Theorem 1.3]{TV2}
(we refer to that work for the definition of a
multi-fold):
\begin{theorem}
\label{orbthm3}
Let $(M_i, g_i)$ be a sequence of unit-volume metrics of class 
(a), (b) or (c) on smooth, closed four-dimensional manifolds $M_i$ satisfying
\begin{align}
\int_{M_i} |Rm_{g_i}|^2 dV_{g_i} \leq \Lambda,
\end{align}
where $\Lambda$ is a constant. 

Assume that
\begin{align} 
\label{a_1}
Vol( B(q,s)) &\geq V_0 s^4, \mbox{ for all } q \in M, 
\mbox{ and } s \leq diam(M)/2,\\
\label{a_3}
b_1(M_i) &< B_1,
\end{align}
where $C_1, \Lambda, B_1$ are constants.
Then a subsequence converges to a 
limit metric space $(M_{\infty}, g_{\infty})$ 
which is a compact,
connected, critical Riemannian multi-fold.
The convergence is smooth away 
from a finite singular set. 

 If we assume instead that
\begin{align} 
\label{a_1'}
& {C_S}_{g_i}(r)  < C_1, \mbox{ for } r < diam(M)/2,
\end{align}
then the same conclusion is true. 
\end{theorem} 

\begin{remark}  The first Betti number assumption
from \cite[Theorem 1.1]{TV} (which allows non-compact 
limits) can similarly be removed, and the Sobolev 
constant assumption relaxed to lower volume growth.
\end{remark}

 Let $R_g$ denote the scalar curvature 
of the metric $g$. Using the above in case (a), we note
the following special corollary.
\begin{corollary}
\label{asdtheorem}
Let $(M, g_i)$ be a sequence of unit volume 
constant scalar curvature anti-self-dual metrics on a fixed closed 
$4$-manifold $M$. Assume that 
\begin{align}
|R_{g_i}| &< C,\\
\label{a_1asd}
Vol( B(q,s)) &\geq V_0 s^4, \mbox{ for all } q \in M, \mbox{ and } s \leq diam(M)/2,
\end{align}
where $C, V_0$ are constants.
Then a subsequence converges to a 
limit space $(M_{\infty}, g_{\infty})$ which is a compact,
connected, anti-self-dual Riemannian multi-fold.
\end{corollary} 

\subsection{Acknowledgements}
The authors would like to thank Gilles Carron, John Lott and 
Joao Santos for enlightening discussions. We also  
thank Edward Fan and Xiuxiong Chen for insightful 
questions and comments on our previous work.

%%%%%%%%%%%%%%%%%%%%%%%%%%%%%%%%%%%%%%%%%%%%%%%%%%%%%%%%%
%%%%%%%%%%%%%%%%%%%%%%%%%%%%%%%%%%%%%%%%%%%%%%%%%%%%%%%%%
\section{Proof of Theorem \ref{bigthm_j}}
%%%%%%%%%%%%%%%%%%%%%%%%%%%%%%%%%%%%%%%%%%%%%%%%%%%%%%%%%
%%%%%%%%%%%%%%%%%%%%%%%%%%%%%%%%%%%%%%%%%%%%%%%%%%%%%
 The space of $L^2$-harmonic $k$-forms 
$\mathcal{H}^k(M)$ is defined to be those  
$\omega \in \Lambda^k(T^*M)$ satisfying 
$\Delta \omega = 0$, and $\omega \in L^2(M)$.
It is well-known that $ \mathcal{H}^k(M) 
\simeq \overline{H}^k_{(2)}(M)$, the {\em reduced} 
$L^2$-cohomology, see \cite{Carron}. 
We next quote the following finiteness theorem. 

\begin{theorem}[Carron, \cite{Carron2}, \cite{Carron}] 
\label{Carronthm}
Let $(M,g)$ be a complete Riemannian manifold satisfing 
for $p > 2$ the Sobolev inequality
\begin{align}
\label{Csobv}
\mu_p(M) \big( \int_M |u|^{\frac{2p}{p-2}}(x) dV_g \Big)^{1 - \frac{2}{p}}
\leq \int_M |du|^2 (x) dV_g,  \mbox{ for all } u \in C^{\infty}_{cpt}(M).
\end{align}
If the negative part of the Ricci satisfies 
\begin{align}
\label{Cric}
\int_M |Ric_-|^{p/2} dV_g < \infty,
\end{align}
then $\mathcal{H}^1(M)$ is finite dimensional.
If the full Riemannian tensor curvature satisfies
\begin{align}
\label{Criem}
\int_M |Rm|^{p/2} dV_g < \infty,
\end{align}
then $\mathcal{H}^k(M)$ is finite dimensional
for each $1 \leq k \leq n$. 
\end{theorem}
\begin{remark} 
For convenience of the reader, we give an indication of the proof
in \cite{Carron}. A crucial estimate is that the Sobolev constant 
bound yields an estimate on the heat kernel: 
there exists a constant $C$ such that 
for all $x \in M$, and any $ t > 0$,
\begin{align}
k(t,x,x) \leq C t^{-n/2}.
\end{align}
It is then shown that $\mathcal{H}^k(M)$ is finite dimensional, 
using the Cwickel-Lieb-Rosenbljum estimate adapted to the Riemannian 
setting. Similar results were also obtained in \cite{BerardBesson}. 
The Weitzenbock formula for a $1$-form is 
\begin{align*}
\nabla^* \nabla  = \Delta_H  + Ric,
\end{align*}
so for the case of harmonic $1$-forms, only the Ricci assumption 
(\ref{Cric}) is necessary. For $k > 1$, the Weitzenbock 
formula depends upon the full curvature tensor, 
which is why the full curvature assumption (\ref{Criem}) 
is required. We note that this method also gives an explicit estimate on 
the dimension of $\mathcal{H}^k(M)$ in terms of the 
Sobolev constant $\mu_p(M)$, and the $L^{p/2}$ curvature 
integral.
\end{remark}

 We recall a definition
\begin{definition}  We say a component $A_0(r_1,r_2)$ of an annulus 
$A(r_1, r_2) = \{ q \in M \ | \ r_1 < d(p, q) < r_2 \}$
is {\em{bad}} if $S(r_1) \cap \overline{A_0(r_1,r_2)}$ has 
more than $1$ component, where $S(r_1)$ is the 
sphere of radius $r_1$ centered at $p$. 
\end{definition}

  For an annulus  $A_0(r_1,r_2)$, we call a component 
of $S(r_1) \cap \overline{A_0(r_1,r_2)}$ an {\em{inner}} 
sphere. Note that this may have several components --
indeed, this is one of the main subtleties in proving 
the volume growth estimate in Theorem \ref{bigthm_j}.  

\begin{theorem}
Under the assumption in Theorem \ref{bigthm_j}, 
there exists a constant $N_0$ (depending upon $C_S$ and $\Lambda$) such that if 
$A_1, A_2, \dots, A_N$ is a collection of disjoint 
connected bad annuli, then $N <  N_0$. 
\end{theorem}
\begin{proof}
From the assumptions in Theorem \ref{bigthm_j},  the 
estimate (\ref{Csobv}) is satisfied for $p = n$,
therefore Theorem \ref{Carronthm} says that 
$\mathcal{H}^1$ is finite dimensional. 
Letting $H^1(M)$ and $H^1_c(M)$ denote the first 
cohomology and first cohomology with compact support, 
respectively, this implies that 
\begin{align}
\label{finim}
\mbox{Image}  \big( H^1_c(M) \rightarrow H^1(M) \big)
\end{align}
has finite dimension, since the space 
in (\ref{finim}) injects into $\mathcal{H}^1$ 
(see \cite{AndersonL2}). The rest of the 
argument just uses the finiteness of (\ref{finim}). 
Let $A_i = A(r_i, s_i)$. Without loss of generality, 
let us assume that the sequence of radii $r_i$ is 
non-decreasing. Since $A_i$ is bad, 
the inner sphere $S(r_i) \cap A_i$ has $N_i > 1$ 
components, 
call them $\mathcal{C}_{i,j}, j = 1 \dots N_i$. 
Take any two components, say $\mathcal{C}_{i,1}$
and $\mathcal{C}_{i,2}$, and
let $p_{i,j}$ be any point of $\mathcal{C}_{i,j}$.
For fixed $i$, let $\gamma_{i,2, 1}(t)$ be a
curve in $A_i$ connecting $p_{i,2}$ with 
$p_{i,1}$. We can always find such a curve 
since, by assumption, $A_i$ is connected.
We can also find a curve $\alpha_{i, 1, 2}(t)$ 
connecting  $p_{i,1}$ with 
$p_{i,2 }$, with image in $B(p, r_i)$.
Joining these curves, we find 
a closed loop $\beta_{i} = \gamma_{i,2, 1} \# \alpha_{i, 1, 2}  
$ based at $p_{i,1}$ 

Next, we find a function $f_{i}$ defined on 
$\overline{A_i}$ such that $f_{i}$ is supported in a
neighborhood of $\mathcal{C}_{i,1}$, with 
$f_i = 1$ on smaller neighborhood of $\mathcal{C}_{i,1}$, 
and $f_i = 0$ in a neighborhood of all other components 
$\mathcal{C}_{i,j}$.
Then the $1$-form $\alpha_i = d f_i$ clearly has
an extension to a smooth $1$-form on $M$, 
which is closed (but not exact), and is supported on $A_i$. 
We claim that 
\begin{align}
\int_{\beta_i} \alpha_i = 1.
\end{align}
To see this, since $\alpha_i$ is supported 
on $A_i$, 
\begin{align}
\int_{\beta_i} \alpha_i = \int_{ \gamma_{i,2, 1}(t) } \alpha_i
=  \int_{ \gamma_{i,2, 1}(t) } df = f(p_{i,1}) - f(p_{i,2}) = 1.
\end{align} 
Furthermore, 
\begin{align}
 \int_{\beta_i} \alpha_j = 0,   \ \ i < j.
\end{align}
This is true since we have assumed the annuli are
indexed by increasing radius, the $\alpha_j$ forms 
are supported either on a different component, 
or on an annulus which is further out. 
We have shown that the $\alpha_i$ define non-zero
independent cohomology classes in 
$\mbox{Image} \big( H^1_c(M) \rightarrow H^1(M) \big)$,
and we therefore have
\begin{align}
N_0 \leq  \mbox{dim}  \big\{ \mbox{Image}\big(H^1_c(M) \rightarrow H^1(M)\big) \big\}.
\end{align}
\end{proof}

\begin{remark} We thank Gilles Carron for providing the 
above argument, which was much simpler than our original 
proof. Note also that
the finiteness of the dimension of (\ref{finim}) will be automatically satisfied 
if either $H^1(M)$, $H^1_c(M)$, or $\mathcal{H}^1(M)$ is 
finite dimensional.
\end{remark}

We quote the following theorem from our previous work.
\begin{theorem}(\cite[Theorem 5.2]{TV2})
\label{bigthm}
Let $(M,g)$ be a complete, non-compact, $n$-dimensional 
Riemannian orbifold (with finitely many singular 
points) with base point $p$.
Assume that there exists a constant $C_1 > 0$ so that
\begin{align}
\label{cond4}
Vol(B(q,s)) \geq C_1 s^n,
\end{align}
for any $q \in M$, and all $s \geq 0$.
Assume furthermore that as $r \rightarrow \infty$,
\begin{align}
\label{decay1}
\underset{S(r)}{sup} \ |Rm_g| &= o(r^{-2}),
\end{align}
where $S(r)$ denotes the sphere of radius $r$ centered
at $p$. If $(M,g)$ contains only finitely many 
disjoint bad annuli, then $(M,g)$ has finitely many ends, 
and there exists a constant $C_2$ so that 
\begin{align}
\label{vga}
Vol(B(p,r)) \leq C_2 r^n,
\end{align}
Furthermore, each end is ALE of order $0$. 
\end{theorem}

Theorem \ref{bigthm_j} follows from the above, since 
we have shown there are only finitely many bad annuli,
and noting that the lower volume growth estimate is implied by the 
Sobolev constant bound (see \cite[Lemma 6.1]{TV}). 
Note also that there is an explicit bound the number of ends  
in terms of $C_S$ and $\Lambda$, see \cite[Theorem 3.3]{Carron2}. 

\begin{theorem}
\label{apporb}
 Theorem \ref{bigthm_j} is valid 
if $(M,g)$ is assumed to be a smooth orbifold 
with finitely many singular points. 
\end{theorem}
\begin{proof}
It can be verified that Carron's arguments are 
valid for smooth orbifolds, and since Theorem \ref{bigthm}
is also valid for orbifolds, the proof is identical to 
the smooth case.
Alternatively, instead of using Carron's result 
for orbifolds, one can argue, albeit non-effectively,
as follows. 
Take a smoothing of $(M,g)$ by 
cutting out a small ball about each orbifold 
singularity and gluing in any smooth metric. 
This can be done because each boundary $S^3 / \Gamma$ 
certainly bounds some smooth manifold. 
We now have a smooth manifold $(\tilde{M}, \tilde{g})$ 
which is isometric to $(M,g)$ outside of a
large ball $B(p, R)$. 
The manifold  $(\tilde{M}, \tilde{g})$ has the 
same asymptotic behaviour as the original $(M,g)$,
so all of the assumptions of Theorem \ref{bigthm_j}  
are satisfied by $(\tilde{M}, \tilde{g})$.
Applying Theorem \ref{bigthm_j}, we obtain an upper volume estimate 
for balls in $(\tilde{M}, \tilde{g})$, which clearly implies 
an upper volume estimate for the original $(M,g)$,
since the asymptotics are the same. 
\end{proof}
%%%%%%%%%%%%%%%%%%%%%%%%%%%%%%%%%%%%%%%%%%%%%%%%%%%%%%%%%%%%%%%%

\section{Sobolev constant and local regularity}
 Recall from the introduction that we have the following notion 
of local Sobolev constant. 
For $p \in M$, and $r > 0$, we define $C_{S}(p, r)$ to be the best 
constant such that  
\begin{align}
\label{mainSob3'}
\Vert f \Vert_{L^{4}} \leq C_{S}(p,r)   
\Vert \nabla f \Vert_{L^2},
\end{align}
for all $f \in C^{0,1}$ with compact support 
in $B(p, r)$. Clearly 
\begin{align}
\label{sob+}
C_S (p, r) \mbox{ is an increasing 
function of } r,
\end{align}
and
\begin{align}
\label{soblim}
\lim_{r \rightarrow 0} C_S(p,r) = C_E,
\end{align}
where $C_E$ is the best constant for the 
Euclidean Sobolev inequality.
\begin{proposition}
\label{sobregthm}
Assume that (\ref{generaleqn}) is satisfied, 
and let $B(p_0,2r)$ be a geodesic 
ball around the point $p_0$.
Assume there exists a constant $V_0 > 0$ so 
that 
\begin{align}
\label{asp1}
Vol( B(q,s)) \geq V_0 s^4
\end{align}
for all $q \in B(p_0, r),$ and $s \leq r$.

Then there exists a 
constant $\epsilon_0$  (depending upon $V_0$) so that if 
\begin{align*}
\Vert Rm \Vert_{L^2(B(p_0,2r))} = 
\left\{ \int_{B(p_0,2r)} |Rm|^2 dV_g \right\}^{1/2} \leq \epsilon_0,
\end{align*}
then 
\begin{align}
C_S(p_0, r/2) \leq C V_0^{-1/4},
\end{align}
and $C$ does not depend on $V_0$.  

\end{proposition}
\begin{remark} We conjecture that for (\ref{generaleqn}), $\epsilon_0$ can moreover be
taken independent of $V_0$ and assumption (\ref{asp1}) is not
necessary. This was recently proved for Einstein metrics 
by Cheeger-Tian \cite{CheegerTiannew}. In fact,
in \cite[Section 11]{CheegerTiannew} this 
was already conjectured to hold for anti-self-dual metrics and 
K\"ahler metrics with constant scalar curvature.   
\end{remark}
\begin{proof}
Without loss of generality, rescale so that $r=1$.
The proof goes by contradiction. Assume we have a sequence 
of critical metrics $g_i$, $i = i \dots \infty$,  
and a sequence $\epsilon_i \rightarrow 0$ as $i \rightarrow \infty$, with 
\begin{align*}
\Vert Rm \Vert_{L^2(B(p_0,2))} \leq \epsilon_i,
\end{align*}
and that 
\begin{align}
C_S(p_0, 1/2) > C V_0^{-1/4}
\end{align}
(we will choose $C$ later). 

We first choose a ``nice'' ball, one for which 
the Sobolev constant is controlled in the ball, and 
also for nearby points. The following lemma 
is for a fixed metric. 
\begin{lemma}
There exist a point $p_{\infty} \in B(p,1)$
and a radius $0 < r_{\infty} \leq 1/2$ 
such that $C_S(p_{\infty}, r_{\infty}) > C V_0^{-1/4}$, 
and $C_S( p, r_{\infty}/2) \leq C V_0^{-1/4}$ for all $p \in 
B(p_{\infty}, r_{\infty})$.
\end{lemma}
\begin{proof}
From assumption, we have $C_S(p_0, 1/2) > C V_0^{-1/4}$.
If $C_S(p, 1/4) \leq C V_0^{-1/4}$ for all $p \in B(p_0, 1/2)$, 
then let $B = B(p_0, 1/2)$. 
Otherwise, there exists a point $p_1 \in B(p_0, 1/2)$ 
with $C_S(p_1, 1/4) > C V_0^{-1/4}$. 
If $C_S(p, 1/8) \leq C V_0^{-1/4}$, for all $p \in B( p_1, 1/4)$, 
then we let $B = B(p_1, 1/4)$.
We continue inductively, 
assume we have chosen $p_{i-1}$ 
with $p_{i-1} \in B(p_{i-2}, 2^{-i+1})$ and 
$C_S( p_{i-1}, 2^{-i}) > C V_0^{-1/4}$.
If $C_S(p, 2^{-i}) \leq C V_0^{-1/4}$ for all 
$p \in B( p_{i-1}, 2^{-i})$ then we stop, 
and let $B =  B(p_{i-1},2^{-i})$.
Otherwise, there exists a
point  $p_i \in B(p_{i-1}, 2^{-i})$ and 
$C_S ( p_i, 2^{-i-1}) > C V_0^{-1/4}$. 

 We claim this procedure must stop in finitely 
many steps. We have
\begin{align*}
d(p_0, p_i) &\leq d(p_0, p_1) + \sum_{j=1}^{i-1} d(p_i, p_{i+1}) \\
& < 1/2  + \sum_{j=1}^{i-1} 2^{-j-1}\\
& < 1/2 + 1/2 = 1. 
\end{align*}
The sequence of points $\{ p_i \}$ are therefore all contained 
in the unit ball $B(p,1)$. If the process did not 
stop,  then there would exists a limit point $q$.
Our first restriction on $C$ is that 
$ C V_0^{-1/4} > C_E$. But from the the choices, we would 
clearly have $\lim_{r \rightarrow 0} C_S(q,r) > C_E$,
which contradicts (\ref{soblim}).  
\end{proof}

Now we apply the lemma to each metric $g_i$, to find 
points $p_i \in B(p,1)$ (ball in the $g_i$ metric)
and $r_i < 1/2$, with
$C_S(p_i, r_i) > C V_0^{-1/4}$, and 
$C_S(p, r_i/ 2) \leq C V_0^{-1/4}$ for all $p \in 
B(p_i,r_i)$ (Sobolev constant with respect to the 
$g_i$ metric). 
 
 Now we rescale to make $r_i$ unit size, that 
is, define $\tilde{g}_i = r_i^{-2} g_i$. 
and consider the sequence of balls 
$\tilde{B}(p_i,2)$. By scale invariance, we have 
\begin{align*}
\Vert Rm \Vert_{L^2(\tilde{B}(p_i,2))} \leq \epsilon_i.
\end{align*}

By our previous $\epsilon$-regularity theorem
\cite[Theorem 3.1]{TV}, 
the curvature and all of its derivatives are uniformly 
bounded in $B(p_i, 3/2)$
\begin{align}
|\nabla^k Rm| \leq C_k \epsilon_i,
\end{align}
where $C_k$ depends upon $k$ and $C V_0^{-1/4}$. 
Since $\epsilon_i \rightarrow 0$, we use 
the theorem of Cheeger-Gromov to extract 
a flat limit space, $B(p_{\infty}, 3/2)$,
with smooth convergence to the limit. From the assumption 
$Vol(B(p, r)) \geq V_0 r^n$, and since the limit space
is flat, we therefore get a bound on the Sobolev constant 
of the limit space, $C_S \leq C_1 V_0^{-1/4}$. 
Choosing $C > C_1$, we arrive at a 
contradiction. 

\end{proof}
Theorem \ref{higherlocalregthm} then follows from 
Proposition \ref{sobregthm} and \cite[Theorem 3.1]{TV}.

%%%%%%%%%%%%%%%%%%%%%%%%%%%%%%%%%%%%%%%%%%%%%%%%%%%%%%%%%%%
%%%%%%%%%%%%%%%%%%%%%%%%%%%%%%%%%%%%%%%%%%%%%%%%%%%%%%%%%%%
\section{Main volume estimate}

 In this section, we discuss the proof of Theorems
\ref{orbthm2} and \ref{orbthm3}, and 
Corollary \ref{asdtheorem}. We first prove Theorem \ref{orbthm2}. 

\begin{proof}(of Theorem \ref{orbthm2}). 
The proof is based on the argument from \cite{TV2}, 
with some modifications. First, let us assume that the volume growth 
estimate from Theorem \ref{orbthm2}  holds 
for $r \leq r_0$, where $r_0$ is some fixed scale. 
That is, let us assume that 
\begin{align}
\label{volloc}
Vol(B_{g_i}(p,r)) \leq V r^4
\end{align}
 for all
$ p \in M_i$, and all $r \leq r_0$. 

 From the $\epsilon$-regularity Theorem \ref{higherlocalregthm},
$(M_j, g_j, p_j)$ converges to a limiting multi-fold 
$(M_{\infty}, g_{\infty}, p_{\infty})$ for some subsequence 
$\{j \} \subset \{i \}$, with finitely many 
$C^0$-multi-fold singular points (recall that a $C^0$-multi-fold 
point means that for each cone at a singular point, 
the metric has a continuous extension to the universal cover of the 
punctured cone). The argument for this is the 
same as in \cite{TV2}. 
\begin{proposition}
The singular points are smooth orbifold points. 
That is, if $x$ is a singular point, then 
for some $\delta >0$, the universal cover of 
$B(x, \delta) \setminus \{x\}$ is diffeomorphic to a punctured 
ball $B^4 \setminus \{0\}$ in $\mathbb{R}^4$, and the lift of $g$, 
after diffeomorphism, extends to a smooth critical metric $\tilde{g}$
on $B^4$.
\end{proposition}
\begin{proof}
In the case we assume a bound on the Sobolev 
constant (but no Betti number bound), this 
follows directly from \cite[Theorem 6.4]{TV2}.
In the case where we only assume a lower volume 
growth bound, we claim that for any singular 
point $p$, there exists a constant $C_s$ so that
\begin{align}
 \Vert u \Vert_{L^4( B(p,\epsilon))} 
\leq C_s  \Vert \nabla u \Vert_{L^2(  B(p,\epsilon))   }
, \  u \in C^{0,1}_{cpt}(   B(p,\epsilon) ).
\end{align}
Indeed, in a neighborhood of a singular point, 
a $C^0$-orbifold is just a $C^0$-perturbation of a 
flat cone, so clearly it satisfies a Sobolev inequality. 
The result then follows again directly from \cite[Theorem 6.4]{TV2}.
\end{proof}

\begin{proposition}
\label{numprop} 
If (\ref{volloc}) is satisfied, and $(M_j, g_j, p_j)$ 
converges to a limiting multi-fold 
$(M_{\infty}, g_{\infty}, p_{\infty})$, as $j \rightarrow
\infty$, then 
there is a bound on the number of cones at a 
singular point of convergence, depending only upon 
$\Lambda$, $V_0$, and $B_1$ in the first case, and 
depending only upon $\Lambda$ and $C_S(r)$ in the 
second case. The bound does not depend upon the 
constant $V$ in (\ref{volloc}). 
\end{proposition}
\begin{proof}
We consider the first case. 
At any limiting multi-fold point $p_{\infty}'$, for 
$\delta > 0$ small,  we look
at the balls $U_j = B(p_j', \delta) \subset M_j$ where
$p_j' \rightarrow p_{\infty}'$ as $j \rightarrow \infty$. 
The are manifolds with boundary components spherical 
space forms, with the metrics arbitrarily close to 
the limiting multi-fold cone metric. 
We apply the Gauss-Bonnet Theorem to conclude
\begin{align}
\label{GB0}
\chi(U_j) \sim \int_{U_j} Rm * Rm + \sum_{k} \frac{1}{|\Gamma_k|}. 
\end{align}
where $Rm * Rm$ is a quadratic curvature expression, 
the second sum is over the boundary components, and $\Gamma_k$ denotes 
the orbifold group for each cone.
Note the final term is an approximation of the boundary integral, 
with vanishing error term as $\delta \rightarrow 0$. 
In terms of Betti numbers,
\begin{align}
\label{GB}
1 - b_1(U_j) + b_2(U_j) - b_3(U_j)
\sim \int_{U_j} Rm * Rm + \sum_{k} \frac{1}{|\Gamma_k|}.
\end{align}
We write $M_j = U_j \cup V_j$, where 
\begin{align}
U_j \cap V_j  \sim \coprod_k  S^3 / \Gamma_k. 
\end{align}
The Mayer-Vietoris sequence for $U_j$ and $V_j$ is 
\begin{align} 
0= H_2( U_j \cap V_j) \rightarrow H_2( U_j) \oplus H_2 (V_j) 
\rightarrow H_2( M_j) \rightarrow H_1( U_j \cap V_j) = 0, 
\end{align}
which gives $b_2 (M_j) = b_2(U_j) + b_2(V_j)$.

Using (\ref{GB}), we can estimate 
\begin{align}
\sum_{k} \frac{1}{|\Gamma_k|} \leq 1 + b_2(U_j) + C \cdot \Lambda \leq 1 + b_2(M_j) 
+C \cdot \Lambda.
\end{align}
Applying the Gauss-Bonnet formula to the oriented manifold $M_j$, 
\begin{align}
2 - 2 b_1(M_j) + b_2(M_j) = \int_{M_j} Rm * Rm,
\end{align}
yields the estimate
\begin{align}
b_2(M_j) \leq C \cdot \Lambda + 2 b_1(M_j) - 2, 
\end{align}
so we have
\begin{align}
\sum_{k} \frac{1}{|\Gamma_k|} \leq 2 C \cdot \Lambda + 2 B_1 - 1. 
\end{align}
The lower volume growth estimate (\ref{a_1}) clearly implies 
the the orders of the orbifold groups are bounded 
strictly from below, so the proposition follows in this case. 

 In the second case, we use Carron's bound on the number 
of ends $N$ of a complete space $X$ \cite[Theorem 0.4]{Carron2}
\begin{align}
\label{ends}
N \leq  1 + C \cdot  C_S \int_X | Ric_-|^2 dx \leq 1 + C \cdot C_S \cdot \Lambda, 
\end{align} 
which, as mentioned above, is also valid for smooth orbifolds.
Since we have a volume growth bound (\ref{volloc}) 
we can perform a standard bubbling analysis.
This analysis was carried out for Einstein metrics
in \cite{Bando}, \cite{Bando2}, see also \cite{NakajimaConvergence} for a nice 
description of this bubbling process.  The same method works 
in our case, with a few modifications. 

  We recall the main steps
of the bubbling analysis. 
  Let $S$ denote the singular set of convergence. 
Note that in contrast to the Einstein case, a 
point $p \in S$ may actually be a smooth point 
of the limit.  
 For $0 < r_1 < r_2$, we let $D(r_1,r_2)$ denote 
$B(p, r_2) \setminus B(p, r_1)$.
Given a singular point $x \in S$, take a 
sequence $x_j \in (M_j, g_j)$ such that 
$\lim_{ j \rightarrow \infty} x_j = x$
and $B(x_j, \delta)$ converges to 
$B(x, \delta)$ for all $\delta > 0$. 
 We choose a radius $r_{\infty}$  sufficiently 
small and the sequence $x_j$ to satisfy
\begin{align}
\underset{B( x_j, r_{\infty})}{ \mbox{sup}} |Rm_{g_j}|^2 = 
|R_{g_j}|^2( x_j) \rightarrow \infty \mbox{ as } j \rightarrow
\infty,  
\end{align}
and 
\begin{align} \int_{B(x,r_{\infty})} |R_{g_{\infty}}|^2 dV_{g_{\infty}}
  \leq \epsilon_0 /2,
\end{align}
where $\epsilon_0$ is the constant in the $\epsilon$-regularity 
Theorem \ref{higherlocalregthm}.

We next choose $r_0(j)$ so  that 
\begin{align}
 \int_{D(r_0, r_{\infty})} |R_{g_j}|^2 dV_{g_j} = \epsilon_0, 
\end{align}
and again $D_j(r_0, r_{\infty}) = B(x_j, r_{\infty}) \setminus 
B(x_j, r_0)$. 
An important note, which differs from the 
Einstein case, the annulus $D(r_0, r_{\infty})$ 
may have several components. 

 Since the curvature is concentrating at $x$, 
$r_0(j) \rightarrow 0$ as $j \rightarrow \infty$, 
the rescaled sequence 
$(M_j, r_0(j)^{-2} g_j, x_j)$ has a subsequence 
which converges to a complete, non-compact 
multi-fold with finitely many singular points, 
which we denote by 
\begin{align}
M_{\infty, i_1}, \ 1 \leq i_1 \leq 
\# \{S\}.
\end{align}
Note that by assumption, $M_{\infty, i_i}$ 
has bounded Sobolev constant. Since 
\begin{align}
\int_{D(1, \infty)} |Rm|^2 dV_g \leq \epsilon_0, 
\end{align}
there are no singular points outside 
of $B(x,1)$. 

 On the noncompact ends, since we are assuming an 
upper volume growth estimate, the proof of \cite[Theorem 1.3]{TV}
allows us to conclude that the metric is ALE  of order $\tau$ for any 
$\tau < 2$. As in \cite[Proposition 4]{Bando}, we conclude 
that the neck regions (for large $j$) will be 
arbitrarily close to a portion of a flat cones 
$\mathbf{R}^4 / \Gamma$, possibly several cones 
if $M_{i_1}$ has several ends (again in contrast 
with the Einstein case). 
The convergence at a singular point $x_{i_1}$ is 
that the ALE multi-fold $M_{\infty, i_1}$ is 
bubbling off, or scaled down to a point 
in the limit, with each end of $M_{\infty, i_1}$ corresponding 
exactly to a cone at a multi-fold point of the limit 
$(M_{\infty}, g_{\infty})$. 
An important fact is that the 
fundamental group of an end and the 
group of the orbifold cone onto which the 
end is glued must be isomorphic (together 
with their actions on $\mathbb{R}^4$), 
this also follows from the proof in \cite{Bando},
see also \cite[Theorem 2.5]{NakajimaConvergence}. 
In particular, the number of ends of $M_{\infty, i_1}$ is the 
same as the number of components of 
$\partial B (x_{i_1}, r_{\infty})$.

 To further analyze the degeneration at the 
singular points, we look at the multi-fold $M_{\infty, i_1}$
with singular set $S_{i_1}$. If $S_{i_1}$ is
empty, then we can stop, as the number of cones
must be bounded. If not, we do the 
same process as above for each singular 
point of $M_{\infty, i_1}$, and obtain ALE multi-folds 
\begin{align}
M_{\infty, i_1, i_2}, \ 1 \leq i_2 \leq \# \{S_{i_1} \},
\end{align}
each of which are just limits of re-scalings of the 
original sequence around the appropriate basepoints.  
If $M_{\infty, i_1, i_2}$ has singularities, then 
we repeat the procedure. 
This process 
must terminate in finite steps, since in this 
construction, each singularity takes at least
$\epsilon_0$ of curvature: 
the $L^2$-curvature bound (\ref{l2a}) clearly implies the number 
of steps in the procedure is bounded by the number $b = \Lambda / \epsilon_0$. 
As pointed out in 
\cite{Bando2}, there could be some 
overlap if any singular point lies 
on the boundary of $B(1)$ at some stage 
in the above construction. But there can 
only be finitely many, and then there 
must also be a singular point in the 
interior of $B(1)$, so we still take 
away at least $\epsilon_0$ of curvature 
at each step. 

  Note that in each step of the bubbling process, 
each end of the multi-fold obtained in the 
$k$th step will be glued to one of the 
cones at a  multi-fold singular point 
of the $(k-1)$st step, along a neck 
region which is close to a portion of a flat 
cone. Again, we use the fact that the 
fundamental group of an end and the  
group of the orbifold cone onto which the 
end is glued must be isomorphic (with 
isomorphic actions on $\mathbb{R}^4$). 
Consequently, the bound $N$ on the number of ends (\ref{ends})
and the bound $b$ on the number of steps in the 
process, together imply that the number of cones at any 
singular point must be bounded by $N^b$.
\end{proof}

\begin{remark} 
As was discussed in \cite[page 369 ]{TV2}, in the 
K\"ahler case only irreducible singular points can occur
in limit. We remark that this still holds under the weaker assumptions in 
this paper. That is, in the K\"ahler case there are never multiple cones at a
singular point of convergence. 
\end{remark}

\begin{remark}The above proposition was proved using a
alternative method in \cite[Proposition 7.2]{TV2}. 
However, that argument required both a Sobolev constant and a 
first Betti number bound. 
\end{remark}

 Next, we have
\begin{proposition} 
\label{evol}
Let $(M,g)$ be a smooth multi-fold with 
finitely many singular points and $g$ a critical 
metric. Assume that the number of singular 
points is uniformly bounded by the number $N_1$,
that the number of cones at any 
singular point is uniformly bounded by 
the number $N_2$, and that there exists a constant $V_0 > 0$ so 
that 
\begin{align*}
Vol( B(q,s)) \geq V_0 s^4
\end{align*}
for all $q \in B(p, 2),$ and $s \leq 2$.
If 
\begin{align*}
\Vert Rm \Vert_{L^2(B(p,2))} = 
\left\{ \int_{B(p,2)} |Rm|^2 dV_g \right\}^{1/2} \leq \epsilon_0,
\end{align*}
then there exists a constant $A_0$ such that
\begin{align}
Vol(B(p,1)) \leq A_0,
\end{align}
where $A_0$ depends only upon $N_1$, $N_2$ and $V_0$. 
\end{proposition}
\begin{proof}
If $(M,g)$ is smooth, then by  Theorem \ref{higherlocalregthm},
\begin{align}
\label{uest}
\underset{B(p, 1)}{sup}| Rm| \leq \frac{1}{4} C \epsilon_0. 
\end{align}
By Bishop's volume comparison theorem, we must have
$Vol(B(p,1)) \leq A'$, where $A'$ depends only
lower volume growth constant $V_0$ (since $\epsilon_0$ 
only depends on $V_0$).  

 In the case that $(M,g)$ is a smooth orbifold, 
we claim that the $\epsilon$-regularity 
theorem still holds in this setting. 
This is because the argument in \cite[Theorem 3.1]{TV}
uses integration by parts. One then performs a similar 
argument, by cutting out small balls of radius $\delta$ 
around the singular points and verifing that 
the resulting boundary terms vanish as $\delta \rightarrow 0$
(note that it is crucial that the orbifold points 
be {\em{smooth}} for this to be valid). Furthermore, 
Proposition \ref{sobregthm} is remains valid for a
smooth orbifold with a bounded number of singular points. 
The details, which we omit, are straightforward. 
Consequently,  the estimate (\ref{uest}) holds. 
Next, Bishop's volume comparison theorem 
remains valid for smooth orbifolds (see for example \cite{Borz1}), 
with the same constant (or better) as in the smooth case.  
So for an orbifold (no multiple cone points), we 
still obtain 
\begin{align}
Vol(B(p,1)) \leq A'. 
\end{align}
In the more general situation of a multi-fold, 
since there are at most $N_2$ cones at each of the 
$N_1$ singular points, clearly we obtain the estimate
\begin{align}
Vol(B(p,1)) \leq N_1 N_2 A'\equiv A_0. 
\end{align}
\end{proof}

 We also note the following fact, for any metric, 
\begin{align*}
\lim_{r \rightarrow 0} Vol( B(p,r))r^{-4} = \omega_4,
\end{align*}
where $\omega_4$ is the volume ratio of the Euclidean metric
on $\mathbf{R}^4$. Clearly, $ A_0 \geq \omega_4$.

For any metric $(M,g)$, define the maximal volume ratio as
\begin{align} 
MV(g) = \underset{x \in M, r \in \mathbf{R}^+}{\max}
\frac{ Vol ( B(x,r))}{r^4}.
\end{align}
If the theorem is not true, then there exists 
a sequence of critical manifolds $(M_i,g_i)$, 
with $MV(g_i) \rightarrow \infty$, that is, 
there exist points $x_i \in M_i$, and $t_i \in \mathbf{R}^+$ so that 
\begin{align}
\label{contvol}
Vol (B(x_i,t_i))\cdot t_i^{-4} \rightarrow \infty, 
\end{align}
as $i \rightarrow \infty$. We choose a subsequence (which for 
simplicity we continue to denote by the index $i$) and radii 
$r_i$ so that 
\begin{align}
\label{2wchoice}
2 \cdot A_0 = Vol( B(x_i,r_i)) \cdot r_i^{-4} 
= \underset{r \leq r_i}{ \max} \ Vol( B(x_i,r)) \cdot r^{-4},
\end{align}
We furthermore assume that $x_i$ is chosen so that 
$r_i$ is minimal, that is, the smallest radius for which 
there exists some $p \in M_i$ such that 
$Vol(B_{g_i}(p, r)) \leq 2 A_0 r^4$, for all $r \leq r_i$.

First let us assume that $r_i$ has a subsequence converging to zero.
For this subsequence (which we continue to index by $i$), 
we consider the rescaled metric $\tilde{g}_i = r_i^{-2}g_i$, 
so that $B_{g_i}( x_i, r_i)  = B_{\tilde{g}_i}( x_i, 1)$. 
From the choice of $x_i$ and $r_i$, 
the metrics $\tilde{g}_i$ have bounded volume ratio,
in all balls of unit size.  

From the argument above, some subsequence converges
on compact subsets to a complete 
length space $(M_{\infty}, g_{\infty}, p_{\infty})$ 
with finitely many point singularities.  
The limit could either be compact or non-compact. 
In either case, the arguments above imply that the 
limit is a Riemannian multi-fold. 

\begin{claim} 
\label{clm}
The conclusions of Theorem \ref{bigthm_j} 
hold for the limit $(M_{\infty}, g_{\infty}, p_{\infty})$ 
\end{claim}
\begin{proof}
In the case that $M_{\infty}$ is compact the claim is trivial. 
For $M_{\infty}$ non-compact, the remarks at the end 
of \cite[Section 3]{TV} shows that assumption (\ref{decay1}) is
satisfied. Also, from \cite[Lemma 6.1]{TV}, the Sobolev 
constant bound implies a lower volume growth bound 
(this is also valid for orbifolds), so (\ref{cond4}) is satisfied
in both cases. 

 If we make the assumption on $b_1(M_i)$ (but no Sobolev 
constant assumption), then the proof is
the same as in our previous work \cite[Claim 7.1]{TV2}:
since $b_1(M_i) < \infty$, the limit must have finitely many bad annuli. 
The finiteness of ends and upper volume growth estimate 
then follow from Theorem \ref{bigthm}.
In the case where we assume the Sobolev constant
bound (but no $b_1$ bound), the result is contained in
Theorem \ref{apporb}. 
\end{proof}

 From Claim \ref{clm}, we have that  
$M_{\infty}$ has only finitely many ends, and that there 
exists a constant
$A_1 \geq 2 A_0$ so that 
\begin{align}
\label{poq}
Vol(B_{g_{\infty}}(p_{\infty}, r)) \leq 
A_1 r^4, \mbox{ for all } r > 0.  
\end{align}
  If $M_{\infty}$ is compact, then clearly the estimate 
(\ref{poq}) is valid for some 
constant $A_1 \geq 2 A_0$, since the limit has finite 
diameter and volume, and the estimate holds for $r \leq 1$.    

The inequality 
\begin{align}
\label{ballineq}
\int_{ B_{g_i}( x_i, 2 r_i)} |Rm_i|^2 dV_i > \epsilon_0,
\end{align}
must hold; otherwise, as remarked above, we would have 
$Vol (B_{g_i}( x_i, r_i)) \leq A_0 r_i^4$, which 
violates Proposition \ref{evol}. 

If the $r_i$ are bounded away from zero then 
there exists a radius $t$ so that 
\begin{align}
Vol(B_{g_i}(p,r)) \leq 2A_0 r^4, \mbox{ for all } r \leq t, p \in M_i.
\end{align}
We repeat the argument from the first case, but 
without any rescaling. 
Since the maximal volume ratio is bounded on small scales, 
we can extract an multi-fold limit. The limit can 
either be compact or non-compact, but the inequality 
(\ref{poq}) will also be satisfied for some $A_1$, 
Following the same argument, we find a sequence of balls 
satisfying (\ref{ballineq}). 

We next return to the (sub)sequence $(M_i, g_i)$
and extract another subsequence so that 
\begin{align}
\label{3wchoice}
2600 \cdot A_1 = Vol( B(x_i',r_{i}')) \cdot (r_i')^{-4} 
= \underset{r \leq r_i'}{ \max} \ Vol( B(x_i',r)) \cdot r^{-4}.
\end{align}
Again, we assume that $x_i'$ is chosen so that 
$r_i'$ is minimal, that is, the smallest radius for which 
there exists some $p \in M_i$ such that 
$Vol(B_{g_i}(p, r)) \leq 2600 A_1 r^4$, for all $r \leq r_i$.
Clearly, $r_i < r_i'$.  

Arguing as above, if $r_i' \rightarrow 0$ 
as $ i  \rightarrow \infty$,  then we 
repeat the rescaled limit construction, but now 
with scaled metric $\tilde{g}_i = (r_i')^{-2} g_i$, 
and basepoint $x_i'$. We find a limiting 
multi-fold $(M_{\infty}', g_{\infty}', p_{\infty}')$,
and a constant $A_2 \geq 2600 A_1$ so that 
\begin{align*}
Vol(B_{g_{\infty}'}(p_{\infty}', r)) \leq A_2 r^4 \mbox{ for all } r > 0. 
\end{align*}
For the same reason as above, we must have
\begin{align*}
\int_{ B_{g_j}( x_j', 2 r_j')} |Rm_j|^2 dV_j > \epsilon_0.
\end{align*}
If $r_i'$ is bounded below, we argue similarly, but without any rescaling. 

We next consider the ratio $r_i' / r_i$. 
There are $2$ possible cases. 

Case (i): there exists a subsequence (which we continue to 
index with $i$) satisfying $r_i' < C r_i$ for some constant $C$.

Case (ii) :
\begin{align}
\lim_{i \rightarrow \infty} \frac{r_i'}{r_i} = \infty.
\end{align}

In Case (i) we proceed as follows:
We claim that for $i$ sufficiently large, the 
balls $B(x_i,2 r_i)$ (from the first subsequence)
and $B(x_i',2 r_i')$ (from the second) must be disjoint 
because of the choice in (\ref{3wchoice}). 
To see this, if  $B(x_i,2 r_i) \cap B(x_i', 2 r_i') \neq \emptyset$,
then $  B(x_i', 2 r_i') \subset B(x_i, 6 r_i')$.
Then (\ref{poq}) and (\ref{3wchoice}) yield
\begin{align}
\begin{split}
2600 A_1 (r_i')^4 & = Vol ( B ( x_i', r_i'))\\
& <  Vol ( B ( x_i',2 r_i'))
< Vol (  B ( x_i, 6 r_i')) \leq 2 A_1 ( 6 r_i')^4
 = 2592 A_1 (r_i')^4,
\end{split}
\end{align}
which is a contradiction
(note the last inequality is true for $i$ sufficiently large 
since (\ref{poq}) holds for the limit, which is valid
only in Case (i)).

 In Case (ii) we argue as follows. If the 
balls $B(x_i,2 r_i)$ (from the first subsequence)
and $B(x_i',2 r_i')$ (from the second) are 
disjoint for all $i$ sufficiently large, then we proceed 
to the next step. Otherwise, we look again 
at the scaling so that 
$r_i ' = 1$: $\tilde{g}_i = (r_i')^{-2} g_i$, 
and basepoint $x_i'$.
Then in this rescaled metric, 
\begin{align}
Vol ( B ( x_i', 1)) = 2600 A_1. 
\end{align}
As above, we have a limiting multi-fold 
$(M_{\infty}', g_{\infty}', p_{\infty}')$,
satisfying
\begin{align}
Vol ( B ( p_{\infty}', 1)) = 2600 A_1. 
\end{align}
Proposition \ref{numprop} implies the number 
of cones at a multi-fold point is a priori bounded, 
so from Proposition \ref{evol}, we conclude that
\begin{align}
\label{ballineqaa}
\int_{ B_{g_{\infty}'}(p_{\infty}', 2)} |Rm|^2 dV > \epsilon_0,
\end{align}
There is now a singular point of convergence corresponding 
to the balls $B(x_i, r_i)$ in the first subsequence.
But since we are in Case (ii), in the $g_i'$ metric, these balls 
must limit to a point in $M_{\infty}'$. 
The only possibility is that the original sequence 
satisfied 
\begin{align}
\label{ballineqbb}
\int_{ B_{g_i}(x_i', 2 r_i')} |Rm_{g_i}|^2 dV_i > 2 \epsilon_0,
\end{align}
for all $i$ sufficiently large. 

We repeat the above procedure, considering possible 
Cases (i) and (ii) at each step. At the $k$th step,  
we can always account for at least $k \cdot \epsilon_0$ of 
$L^2$-curvature.  The process must terminate in 
finitely many steps from the bound 
$ \Vert Rm_i \Vert_{L^2} < \Lambda$. This contradicts 
(\ref{contvol}), which finishes the proof.

\end{proof}

\begin{remark} 
In the proof of \cite[Theorem 1.2]{TV2},
we neglected to consider the possibility of Case (ii),
the above fixes this oversight. Another point is
that Propositions \ref{numprop} and 
\ref{evol} depend crucially on the limiting multi-folds being 
{\em{smooth}},  which requires the removable singularity 
result \cite[Theorem 6.4]{TV2}. 
The argument given by Anderson in \cite{Andersonc} misses 
this important point. 
It will be interesting to find a valid proof not 
using the removable singularity theorem.
Note also that in \cite{Andersonc}, Anderson claims to prove 
the upper volume growth estimate in Theorem \ref{orbthm3} 
without requiring either a Betti 
number assumption or a Sobolev constant bound (only assuming
a lower volume growth bound). 
We point out that his argument is incomplete --
in that work insufficient consideration is given to the 
connectedness properties of geodesic spheres and annuli.
Proper consideration of these connectedness properties 
is an absolutely crucial point, as can be observed in our proof.
\end{remark}

  Theorem \ref{orbthm3} is proved in a similar manner as the corresponding 
theorem in \cite{TV2}, using Theorems \ref{bigthm_j}, 
\ref{higherlocalregthm}, and the volume growth estimate in 
Theorem \ref{orbthm2}.

 We next prove Corollary \ref{asdtheorem}, which is
a simple consequence of the Gauss-Bonnet Theorem and 
Hirzebruch Signature Theorem in dimension four (see \cite{Besse}):
\begin{align}
\label{Euler}
8 \pi^2 \chi(M) &= \frac{1}{6}\int_M R^2 - \frac{1}{2} \int_M |Ric|^2 + 
\int_M |W|^2,\\
\label{Signature}
12 \pi^2 \tau(M) &=  \int_M |W^+|^2 - \int_M |W^-|^2.
\end{align}
 In the anti-self-dual case, $W^+ \equiv 0$,  so we have 
\begin{align}
\label{asd1}
8 \pi^2 \chi(M) &= \frac{1}{6}\int_M R^2 - \frac{1}{2}\int_M |Ric|^2 + \int_M |W^-|^2,\\
\label{asd2}
12 \pi^2 \tau(M) &=  - \int_M |W^-|^2.
\end{align}
Add these equations together to obtain 
\begin{align}
\label{added}
8 \pi^2 \chi(M) +12 \pi^2 \tau(M)  &= \frac{R^2}{6}Vol(M) - \frac{1}{2}\int_M |Ric|^2\\
 &= \frac{R^2}{6} - \frac{1}{2}\int_M |Ric|^2,
\end{align}
since the scalar curvature is constant, and $Vol = 1$. 
The topology of the manifold $M$ is fixed,
and the scalar curvature is uniformly bounded,  
so we find that $ \frac{1}{2}\int_M |Ric|^2$ is bounded. 
Also (\ref{asd2}) yields a bound on  $\int_M |W^-|^2$,
so we have the estimate  
\begin{align}
\int_M |Rm|^2 dV_g < \Lambda,
\end{align}
for some constant $\Lambda$. 
Therefore the assumptions of Theorem \ref{orbthm3} are satisfied, 
which finishes the proof. 

%%%%%%%%%%%%%%%%%%%%%%%%%%%%%%%%%%%%%%%%%%%%%%%%%%%%%%%%%%%%%%%%%%%%%
\section{ALE metrics and removable singularities}
%%%%%%%%%%%%%%%%%%%%%%%%%%%%%%%%%%%%%%%%%%%%%%%%%%%%%%%%%%%%%%%%
 A related problem is to find geometric conditions 
so that each end of a complete space will be
ALE of order $\tau > 0$, and to determine the 
optimal order of decay. In \cite{TV} we examined this problem 
for the following cases:\\
\begin{tabular}{ll}\vspace{-3mm}
\\
a. & Self-dual or anti-self-dual metrics with zero scalar curvature.\\

b. & Scalar-flat metrics with harmonic curvature.\\
\end{tabular}\\
\\
 By using Theorems \ref{bigthm_j}, \ref{higherlocalregthm}, 
and the volume growth estimate in Theorem \ref{orbthm2}, 
we have the ollowing improvement  of \cite[Theorem 1.3]{TV}:
\begin{theorem}
\label{decayale}
Let $(M,g)$ be a complete, noncompact 4-dimensional Riemannian 
manifold with $g$ of class (a) or (b)
satisfying, 
\begin{align}
 \int_M |Rm_g|^2 dV_g < \infty.
\end{align}
Assume that
\begin{align}
\label{c_1}
Vol( B(q,s)) &\geq V_0 s^4, \mbox{ for all } q \in M \mbox{ and } s > 0,\\
\label{c_3}
b_1(M) &< \infty,
\end{align}
then $(M,g)$ has finitely many ends, and 
each end is ALE of order $\tau$ for any $\tau < 2$.
If we assume instead that 
\begin{align}
\label{c_1bb}
& C_S < \infty,
\end{align}
then the same conclusion holds. 
\end{theorem}

 To conclude, we mention that the following removable singularity 
theorem for critical metrics is expected:\\

{\it{
\label{remsing}
Let $(M,g)$ be a $C^0$-orbifold with singular point at $x$,
and $g$ be a critical metric satisfying (\ref{generaleqn}). 
Suppose that 
\begin{align}
& \int_{B(x,1)} |Rm_g|^2 dV_g < \infty.
\end{align}
Then the metric $g$ extends to $B(x,1)$
as a smooth orbifold metric. 
That is, for some small $\delta >0$, 
universal cover of $B(x, \delta) \setminus \{x\}$ 
is diffeomorphic to a punctured ball $B^4 \setminus \{0\}$ 
in $\mathbb{R}^4$, and the lift of $g$, after diffeomorphism, 
extends to a smooth critical metric $\tilde{g}$
on $B^4$. }}\\

We plan to address this in a forthcoming paper. 
Such a removable singularity theorem was proved for the special
cases of (a),(b), and (c) in \cite[Theorem 6.4]{TV2}. 
The above generalization is crucial in extending Theorems \ref{orbthm2}, 
\ref{orbthm3}, and \ref{decayale} to the more general 
class of critical metrics satisfying (\ref{generaleqn}),
in particular, Bach-flat metrics. 

\bibliography{VolumeGrowth_references}
\end{document}